\documentclass{amsart}[11pt]

\usepackage[english]{babel}
\usepackage{amsmath,amssymb,enumerate,bbm}
\usepackage[utf8]{inputenc}

\newcommand{\mathd}{\mathrm{d}}

\newcommand{\tmtextit}[1]{{\itshape{#1}}}

\newtheorem{definition}{Definition}
\newtheorem{theorem}{Theorem}
\newtheorem{corollary}{Corollary}
\newtheorem{lemma}{Lemma}

\begin{document}

\title[Continuous--time sparse domination]{Continuous--time sparse domination.}

\author{Komla Domelevo}
\address{(K.D.) Institut de Math\'ematiques de Toulouse, Universit\'e Paul Sabatier, 118 route de Narbonne, F-31062 Toulouse Cedex 9, France}
\email{komla.domelevo@math.univ-toulouse.fr}

\author{Stefanie Petermichl}
\address{(S.P.) Institut f\"{u}r Mathematik, Julius Maximilians Universit\"{a}t W\"{u}rzburg, Emil Fischer Stra{\ss}e 40, D-97074 W\"{u}rzburg, Germany}
\email{stefanie.petermichl@mathematik.uni-wuerzburg.de}

\thanks{Research supported by the ERC project CHRiSHarMa DLV-862402 and by the Alexander von Humboldt Stiftung}

\begin{abstract}
  We develop the self similarity argument known as sparse domination in an abstract martingale setting, 
  using a continuous time parameter. 
  With this method, we prove a sharp weighted $L^p$ estimate for the maximal operator $Y^{\ast}$ of $Y$ 
  with respect to $X$. Here $Y$ and $X$ are uniformly integrable c{\`a}dl{\`a}g Hilbert space
  valued martingales and $Y$ differentially subordinate to $X$ via the square
  bracket process.  We also present a second, very simple proof of the 
  special case $Y=X$, exhibiting the optimal weighted $L^p$ estimate previously 
  only known in a homogeneous context. We point out that in this
  generality, notably including processes with jumps, the special case $Y = X$ addresses a 
  question raised  in the late 70s by Bonami--L{\'e}pingle.
  \end{abstract}

\maketitle

\section{Introduction}

We have a filtered probability space with the usual assumptions. Let $X$ and
$Y$ be adapted uniformly integrable c{\`a}dl{\`a}g martingales. Wang's \cite{Wan1995a}
definition of differential subordination has appeared as the correct
continuous in time replacement of the corresponding notion in probability
spaces with discrete filtration by Burkholder \cite{Bur1991a}. We say that $Y$
is differentially subordinate to $X$ if and only if $[X, X]_t - [Y, Y]_t$ is
non-negative and non-decreasing. Here $[\cdot, \cdot]$ denotes the square bracket. 

\

We develop the notion of sparse domination in this abstract context.  
The idea is to construct self similarity through a stopping time procedure where, 
roughly speaking, the events where the iterated stopping times are finite have 
comparably small measure. This gives rise
to a trajectory-wise domination of the maximal function $Y^{\ast}$ by a 
positive form containing $| X |$ sampled at the consecutive stopping times.
The latter domination lends itself well to an estimate in weighted space
through a change of measure and the use of Doob's inequality.

Indeed,  we call an adapted sequence of stopping times $\{T_i\}_{i\geqslant 0}$ a sparse sequence 
   if  $\{T_i\}_{i\geqslant 0}$ is increasing  with nested sets $E_j =
  \{ T^j < \infty \}$ so that
  $$\forall A^j \subset E_j, A^j \in \mathcal{F}_{T^j} {\text{ there holds }}
    \mathbbm{P} (A^j \cap E_{j + 1}) \leqslant \frac{1}{2} \mathbbm{P}
    (A^j).$$

When $Y^*$ denotes the maximal function of $Y$ we prove the following precise sparse domination. 
There exists a sparse sequence of stopping times depending on $X$ and $Y$, so that 
  there holds almost surely
\[  Y^{\ast}(\omega)\leqslant 8 \sum_{j = 0}^{\infty} | X |_{T^j} (\omega)
    \chi_{E_j} (\omega) .\]
In this context we also call $$\mathcal{S} ( X ) (\omega)= \sum_{j = 0}^{\infty} | X |_{T^j} (\omega)
    \chi_{E_j} (\omega)$$ a sparse operator.  If we denote by $X_{\infty}$ a closure of $X$ 
 then we mean by $|X|$ the martingale $|X|_t=\mathbb{E}(  |X_{\infty} | \, |\, \mathcal{F}_t)$ and $|X|_{T^j}$ is $|X|$ sampled at $T^j$.

Notice that this domination is very powerful in the following sense. 
It holds almost surely and it dominates an object via an easily 
estimated positive quantity, in a context where we 
usually see a high degree of cancellation. 

\

One of the uses of sparse operators are weighted estimates. A weight $w$ is a 
positive uniformly integrable martingale. The weight is in the class $A_p$ if the
quantity $Q_p (w)$, the $A_p$ characteristic of the weight $w$, is finite:
\[ Q_p (w) = \sup_{\tau} \left\| \mathbbm{E} \left[ \left( \frac{w_{\tau}}{w}
   \right)^{\frac1 { p - 1}} \middle| \, \mathcal{F}_{\tau} \right]^{p - 1}
   \right\|_{\infty} , \]
where the supremum runs over adapted stopping times  $\tau$. 

\
We prove the estimate
\[ \| \mathcal{S}(X) \|_{L^p (w)} \lesssim  Q_p (w)^{\max \{1,\frac1{p-1} \}} \| X
   \|_{L^p (w)} .\]

With the help of the sparse domination result,
we deduce the following sharp estimate
\[ \| Y^{\ast} \|_{L^p (w)} \lesssim  Q_p (w)^{\max \{1,\frac1{p-1} \}} \| X
   \|_{L^p (w)}. \]

\

The concept of differential subordination is interesting in its own right.
Its ties to harmonic analysis have a long history and have proven influential, 
especially with ambitious goals such as very precise, sharp norm estimates. 
 Certain classical operators in harmonic analysis,
such as Riesz transforms can be written as a conditional expectation of
certain martingale transforms {\cite{GunVar1979}}. Other, deep connections 
of a probabilistic flavour have surfaced in the last ten to twenty years
{\cite{Pet2000b}}{\cite{Hyt2012a}}. 

\

During the last two decades, much interest has shifted towards precise, sharp weighted 
norm estimates in terms of the characteristic of the
weight. In part this interest begun thanks to the solution of a long standing regularity problem
in PDE through a sharp weighted norm estimate of the Beurling-Ahlfors operator
{\cite{PetVol2002a}}. The study of weights has brought, in a joint effort,  new 
understanding of operators central to harmonic analysis, linking them further to 
operators that act directly on a dyadic grid. Making the passage from dyadic to general filtrations, 
we end up back in probability theory.

\

One of the first sharp weighted norm estimates, was on predictable multipliers for
dyadic filtrations in the interval $[0, 1]$ endowed with Lebesgue measure
{\cite{Wit2000a}}. At the time, this problem had been studied by harmonic analysts using 
the language of wavelet bases. The estimate relied on Bellman functions, together 
with the deep reduction and ideas developed in {\cite{NazTreVol1999a}}. 
Meanwhile, a number of other, beautiful proofs have
appeared, for example {\cite{LacPetReg2010a}}, \cite{CruMarPer2012a}, {\cite{Lac2015a}}. 

\

The first extension of the result to predictable multipliers in general discretely filtered spaces is
subtle {\cite{ThiTreVol2015a}} using a combination of outer
measure space technique and some `small' Bellman functions. 
The solution of the general case, where filtered spaces have a continuous parameter and the
martingales are c{\`a}dl{\`a}g and merely in relation of differential
subordination is by the authors in {\cite{DomPet2016a}}. The proof is technical and has required a `large'
Bellman function with subtle convexity properties, solving the entire problem
at once.

\

These two results estimate $Y$, not $Y^{\ast}$, the maximal
function of $Y$. The first result of this stronger nature is short and elegant: the
sharp weighted norm estimate of the maximal function of predictable
multipliers in the discrete in time filtration case is due to Lacey
{\cite{Lac2015a}} via the new idea of point-wise sparse domination 
with his process from the `top down'.
Lacey's motivation was the (more difficult) domination of Calder{\'o}n-Zygmund 
operators. It appeared in parallel with a similar domination result 
by Lerner-Nazarov {\cite{LerNaz}}. Together, these ideas have caused much 
movement in harmonic analysis. To achieve the domination of Calder{\'o}n-Zygmund operators, 
the language used was the one of dyadic cubes, the atoms in dyadic filtrations. Exponentially 
many (in terms of dimension) dyadic systems were required to achieve the domination. 
It is the main task of this paper to develop a sparse 
domination principle in abstract filtered spaces with continuous time. While it is interesting in its won right, 
we expect numerous applications via a trajectory-wise domination using Brownian filtrations for Riesz transforms. 
The preprint \cite{DahDomPetSkr2019} obtains via this new technique dimension-free estimates (even weighted estimates) 
for Riesz transforms in great generality. We expect numerous other applications of this technique.

\

 Let us summarise the main results of the present paper. We develop sparse 
 operators in an abstract martingale setting. We obtain a trajectory-wise sparse domination result 
  of the maximal function $Y^*$ of $Y$ via the martingale $X$, where $Y$ 
 is differentially subordinate to $X$. The result is very general, there are no restrictions on the 
 continuity of the path of the martingales. Further, the pair of martingales under differential
subordination is in general filtered spaces, thus removing the restrictions of a
discretely filtered space as well as the predictable multiplier property. 
We also pay attention to the special case $Y = X$. It becomes a weighted estimate for the maximal function
of $X$. The first such estimate for the dyadic filtration was proved by
Muckenhoupt {\cite{Muc1972a}}. It was brought to the probabilistic context by
Izumisawa and Kazamaki {\cite{IzuKaz1977a}} for martingales with continuous path. 
Later, Buckley {\cite{Buc1993a}} improved the dependence on the
weight's charateristic in the dyadic setting and Osekowski {\cite{Ose2016a}}
improved the estimate in various ways and provided the proof for martingales with continuous path.

We remark that these results either assume a dyadic filtration with a doubling measure 
or, in continuous time, there is a continuity assumption on the path of the martingales. This is a notable
restriction. There are examples for classical weighted estimates that hold in the case of dyadic filtrations with doubling 
measure but do not hold at all in the general case or with a worse dependence on the characteristic, see \cite{DIPTV2019}. 
In the weighted case, useful facts no longer hold when these
assumptions are dropped. One of these is the self improvement property of the
$A_p$ classes, which is used in some classical proofs.  The openness property states 
that any weight in $A_p$ automatically self improves and is in some better $A_{p-\varepsilon}$. It is an old
observation that this fact is false for weights in general filtrations
allowing jumps {\cite{BekBon1978a}}. An additional homogeneity type
condition controlling the jumps was assumed by Dol\'{e}ans-Dade--Meyer
{\cite{DolMey1979a}} to confirm a weighted maximal inequality.
Bonami--L{\'e}pingle {\cite{BonLep1979a}} observed that for the example weight
that belongs to $A_p$ but no $A_{p - \varepsilon }$ for any $\varepsilon >0$ (i.e. self improvement
fails), the weighted maximal inequality still holds. They phrased the question
for the general filtrations allowing jumps. The positive answer follows for
all $1 < p < \infty$ as a special case of our main result. It is however only
sharp in dependence on the $A_p$ characteristic when $p \leqslant 2$. For this
reason, we give another proof for the special case $Y = X$ that gives the best
estimate for all $p$. Indeed we prove that
\[ \| X^{\ast} \|_{L^p (w)} \lesssim c_p Q_p (w)^{\frac1{p-1}} \| X \|_{L^p
   (w)} . \]

As explained above, our proof yielding the estimate for $Y^{\ast}$ uses sparse
domination. It reduces the weighted estimate to a use of Doob's inequality.

Our proof for the estimate of $X^{\ast}$ uses a very simple domination of the
maximal function, direct, without the use of stopping times, followed by uses of Doob's inequality. This approach is
a modification of a trick due to Lerner {\cite{Ler2008a}}, where it was used
in a different context. 

\

Notice that maximal inequalities involving $X^{\ast}$ can be deduced from the
discrete in time filtered general case through the use of Doob's sampling
theorem, but not the estimates involving $Y^{\ast}$, since differential
subordination is not preserved when sampling the martingale.

\section{Definitions and Main Results}

Let $(\Omega, \mathcal{F}, \mathfrak{F}, \mathbbm{P})$ be a complete filtered
probability space with $\mathfrak{F}= (\mathcal{F}_t)_{t \geqslant 0}$ a
filtration that is right continuous, where $\mathcal{F}_0$ contains all
$\mathcal{F}$ null sets. Let $X$ and $Y$ be uniformly integrable adapted
c{\`a}dl{\`a}g martingales with values in a separable Hilbert space that are
in a relation of differential subordination according to Wang:

\begin{definition}
  $Y$ is called differentially subordinate to $X$ if $[X, X]_t - [Y, Y]_t$ is
  non-negative and non-decreasing in $t$. In this case we also call the ordered pair $(X,Y)$ differentially subordinate.
\end{definition}

For the definition of the square bracket process and its properties, see for example
Dellacherie--Meyer {\cite{DelMey1982a}} or Protter {\cite{Pro2005a}}. Notice
that in particular $[X, X]_0 = | X_0 |^2$ so that differential subordination
of $Y$ with respect to $X$ implies $| Y_0 |^2 \leqslant | X_0 |^2$. 
Recall that for any stopping time $\tau$ the stopping sigma algebra is 
$$\mathcal{F}_{\tau} = \{ \Lambda
\in \mathcal{F}: \Lambda \cap \{ \tau \leqslant t \} \in \mathcal{F}_t \}.$$

\

We make the following definition:

\begin{definition}
   An increasing sequence   $\{ T^j\}_{j\geqslant 0}$ of stopping times  with nested sets $E_j =
  \{ T^j < \infty \}$ is called sparse if 
$$\forall A^j \subset E_j, A^j \in \mathcal{F}_{T^j} {\text{ there holds }}
    \mathbbm{P} (A^j \cap E_{j + 1}) \leqslant \frac{1}{2} \mathbbm{P}
    (A^j).$$
\end{definition}

\begin{definition}
The maximal function associated with $X$ is  $$X^{\ast} = \sup_{t \geqslant 0} | X_t |.$$ 
\end{definition}

We shall throughout this text denote by
the same letters also the closures of the martingales that arise.

\

Here are our main theorems.

\begin{theorem}
\label{Theorem_sparsedomination}
$Y$ differentially subordinate to $X$ then there exists a sparse selection $$(X,Y) \mapsto \{ T^j\}_{j\geqslant 0}$$ such that almost surely
$$Y^*(\omega)\leqslant 8\sum_{j = 0}^{\infty} | X |_{T^j} (\omega)
    \chi_{E_j} (\omega),$$
    where the right hand side has finitely many terms almost surely.
\end{theorem}

The arising sum on the right is a well defined object, also applied to objects where it loses its 
domination property stated in Theorem \ref{Theorem_sparsedomination}. To be precise, there holds $\mathbb{P}(E_j)\to 0$ as $n\to\infty$.

\begin{definition}
If $\{ T^j\}_{j\geqslant 0}$ is a sparse sequence and $E_j = \{ T^j < \infty \}$ its associated sequence of nested sets, we call
$$\mathcal{S}: X\mapsto \mathcal{S} ( X ) = \sum_{j = 0}^{\infty} | X |_{T^j} 
    \chi_{E_j} $$ a sparse operator.
\end{definition}

With the help of the sparse domination in Theorem \ref{Theorem_sparsedomination} we prove a weighted maximal inequality. Let us recall the definition of the $A_p$ class.

\begin{definition}
We call a positive uniformly integrable martingale $w$ a weight. The
quantity $Q_p (w)$ below is the $A_p$ characteristic of the weight $w$. If $Q_p (w)$ finite, then we say $w \in A_p$.
\[ Q_p (w) = \sup_{\tau} \left\| \mathbbm{E} \left[ \left( \frac{w_{\tau}}{w}
   \right)^{\frac1{p-1}} \middle| \, \mathcal{F}_{\tau} \right]^{p - 1}
   \right\|_{\infty} = \sup_{\tau} \left\| w_{\tau} u^{p - 1}_{\tau} \right\|_{\infty} ,\]
where the supremum runs over adapted stopping times  $\tau$ and where 
we write $u^p w = u$.
\end{definition}

First, we prove a weighted estimate for the sparse operator:
\begin{theorem}\label{Theorem_SX}
There exists $c_{\ref{Theorem_SX},p}\geqslant 0$ such that for all functions $X\in L^p(w)$ there holds
  \[ \| \mathcal{S}(X) \|_{L^p (w)} \leqslant c_{\ref{Theorem_SX},p} 
    Q_p (w)^{\max \{1, \frac1{p - 1}\}} \| X \|_{L^p (w)} . \]
  The estimate is sharp in terms of the dependence on $Q_p (w)$. 
\end{theorem}

\begin{theorem}\label{Theorem_YstarX}
  There exists $c_{\ref{Theorem_YstarX},p}\geqslant 0$ such that for all pairs $(X,Y)$ where $Y$ differentially subordinate to $X$ there holds
  \[ \| Y^{\ast} \|_{L^p (w)} \leqslant c_{\ref{Theorem_YstarX},p} 
  Q_p (w)^{\max  \{1, \frac1{p - 1}\}} \| X \|_{L^p (w)} . \]
  The estimate is sharp in terms of the dependence on $Q_p (w)$.
\end{theorem}

In the special case $Y = X$ we also prove by a different method

\begin{theorem}
  \label{Theorem_XstarX}
  There exists $c_{\ref{Theorem_XstarX},p}\geqslant 0$ such that for all martingales $X$ there holds
   \[ \| X^{\ast} \|_{L^p (w)} \leqslant c_{\ref{Theorem_XstarX},p} Q_p (w)^{  \frac1{p - 1}} \| X \|_{L^p (w)} . \]
 The estimate is sharp in terms of the dependence on $Q_p (w)$. 
\end{theorem}

Writing $\frac1{p}+\frac1{p'}=1$, we have in the theorems above, 
$$c_{\ref{Theorem_SX},p}=
\left\{
\begin{array}{ll}
8\sqrt{2}\,2^{\frac1{p-1}} \left(\frac{p^{p'}}{p-1}\right)^{\frac{p-2}{p-1}}
& 1<p<2\\
8 & p=2\\
8\sqrt{2}\,2 \left(\frac{p'^{p}}{p'-1}\right)^{\frac{p-2}{p-1}}
&p>2
\end{array}
\right.,
$$ 

as well as $c_{\ref{Theorem_YstarX},p}=8\,c_{\ref{Theorem_SX},p},$
and  $c_{\ref{Theorem_XstarX},p}=\frac{p^{p'}}{p - 1}$. 
We notice that $c_{\ref{Theorem_SX},p}$ and $c_{\ref{Theorem_YstarX},p}$ are of the form $O(p)$ for large $p$ but explode badly as $p\to 1$.

\

The weighted estimates are well known to be sharp in terms of the dependence on the
$A_p$ characteristic, already for dyadic filtration on $[0, 1]$ endowed with
Lebesgue measure. In this paper we focus on the upper estimates.

\section{Sparse Domination}

In this section we prove Theorem \ref{Theorem_sparsedomination}. Without loss of generality $X$ has non-zero closure and $\|X\|_1>0$.

\

We will use the following preliminary weak type estimate due to Wang
{\cite{Wan1995a}}.  

\begin{theorem}[Wang] $Y$ differentially subordinate to $X$ then for all $\lambda >0$ there holds
  \[ \mathbbm{P} (\{ \omega \in \Omega : (|Y|  + |X|)^{\ast}
     (\omega) > \lambda \}) \leqslant \frac2{\lambda} \| X \|_1 . \]
\end{theorem}
Notice that it trivially implies
\begin{corollary}[Wang]
 $Y$ differentially subordinate to $X$ then  for all $\lambda >0$ there holds
  \[ \mathbbm{P} (\{ \omega \in \Omega : Y^{\ast} (\omega) \vee X^{\ast}
     (\omega) > \lambda \}) \leqslant \frac2{\lambda} \| X \|_1 . \]
\end{corollary}
It also implies the following.
\begin{lemma} \label{Theorem_WangWeakType}
 $Y$ differentially subordinate to $X$ then for all $A\in \mathcal{F}_0$ and for all $\lambda >0$ there holds
  \[ \mathbbm{P} (\{ \omega \in A : Y^{\ast} (\omega) \vee X^{\ast}
     (\omega) > \lambda | X |_0(\omega) \}) \leqslant \frac2{\lambda} \mathbb{P}(A) . \]
\end{lemma}
\begin{proof}
Let us write 
$$\tilde{X} = \chi_{A \cap \{ |X|_0>0\}}X/|X|_0 \text{ and } \tilde{Y} = \chi_{A\cap \{|X|_0>0\}}Y/|X|_0.$$
Notice first that $|X|_0$ is measurable in $\mathcal{F}_0$. So the pair $(\tilde{X},\tilde{Y})$ are martingales under differential subordination. 
Further, notice that $|X|_0(\omega) =0 \Rightarrow |X|_t(\omega) =0 \; \forall t\geqslant 0$. In particular, all future 
increments of $X$ are zero and thus $|X_t(\omega)| \vee |Y_t(\omega)| =0 \; \forall t>0$. Thanks to this and Wang's Theorem we can estimate for all $\lambda >0$
\begin{eqnarray*} 
&&\mathbbm{P} (\{ \omega \in A :
   Y_t^{ \ast} \vee X_t^{ \ast} > \lambda|X|_0 \})
 =  \mathbbm{P} (\{ \omega \in \Omega :
   \tilde{Y}_t^{ \ast} \vee \tilde{X}_t^{ \ast} > \lambda \})\\
& \leqslant & \frac{2}{\lambda} \| \tilde{X} \|_1
= \frac{2}{\lambda}
   \mathbbm{E} \left[ \frac{\chi_{A\cap \{ |X|_0>0\}}}{| X |_0}\mathbbm{E} [| X | \, |\,  \mathcal{F}_0]
   \right] 
= \frac{2}{\lambda} \mathbb{P}(A). 
\end{eqnarray*}

\end{proof}

\subsubsection{Stopping procedure}
Let us select an increasing sequence $\{T^j\}_{j\ge0}$ of stopping times associated to the pair $X$ and $Y$ inductively. 
To do so, start with 
$$T^0(\omega)=\inf\{t>0: |X|_t(\omega)>0\},$$ that is, 
$T^0(\omega)= 0$ if $|X|_0(\omega) >0$ and $\infty$ else. There holds   
$$E_0=\{T^0<\infty \}=\{ \omega \in \Omega: |X|_0(\omega)>0\}.$$
Notice that since $\|X\|_1>0$ we have $\mathbb{P}(E_0)>0$ and  $|X|_{T_0}>0$ on $E_0$. 
Let us define $\mathfrak{F}^0=(\mathcal{F}_t^0)_{t\geqslant 0}$ by $\mathcal{F}^0_t =\mathcal{F}_{T^0 \vee t}$
and consider the martingales 
$$Y^0_t =\chi_{E_0}Y_t  \text{ and } X^0_t = \chi_{E_0}X_t .$$
   
Now we proceed with the iteration. Assume $n\geqslant 1$.
Assume we have filtrations $\mathfrak{F}^0,...,\mathfrak{F}^{n-1}$, an increasing sequence of 
stopping times $T^0,...,T^{n-1}$ with associated nested sets $E_0,...E_{n-1}$ measurable in $\mathcal{F}_0,...,\mathcal{F}_{n-1}$ respectively and pairs of martingales 
$(X^0,Y^0),...,(X^{n-1},Y^{n-1})$ under differential subordination. 

We set the stopping time $T^n (\omega)$ by
\[ T^n (\omega) = \inf \{ t > 0 : Y^{n-1}_t  (\omega)\vee X^{n-1}_t (\omega)>4|X|_{T^{n-1}}(\omega) \}. \]
Notice that $|X|_{T^{n-1}}(\omega) =0 $ implies $|X|_t(\omega) =0$ for all $t\geqslant T^{n-1}$. In particular, all future 
increments of $X$ are zero and thus $|X^{n-1}_t(\omega)| \vee |Y^{n-1}_t(\omega)| =0 \; \forall t>T^{n-1}$. There holds thus $T^n(\omega)= \infty$ if $|X|_{T^{n-1}}(\omega) =0$. Thus
\[ E_n = \{ T^n<\infty \}=\{ \omega \in \Omega : Y^{n-1 \ast}(\omega) \vee X^{n-1 \ast}(\omega) > 4|X|_{T^{n-1}}(\omega) \} .\]
 Let
the filtration $\mathfrak{F}^n = (\mathcal{F}_t^n)_{t \geqslant 0} =
(\mathcal{F}_{T^n \vee t})_{t \geqslant 0}$. Observe that $E_n \in
\mathcal{F}_{T^n}$. We will have to take care of the foot of the next pair of martingales in the presence of a jump just prior to a stopping time. 
By differential subordination almost surely 
$$\left| Y^{n-1}_{T^n} (\omega) - Y^{n-1}_{T^n_-}(\omega) \right| \leqslant \left| X^{n-1}_{T^n} (\omega) - X^{n-1}_{T^n_-} (\omega) \right|.$$ 
Thus, there exists a linear operator $r_{T^n} (\omega) \in \mathcal{F}_{T^n}$ such that $| r_{T^n} |
\leqslant 1$ and 
$$Y^{n-1}_{T^n} (\omega) - Y^{n-1}_{T^n_-} (\omega) = r_{T^n} (\omega)  \left(X^{n-1}_{T^n}(\omega) - X^{n-1}_{T^n_-} (\omega)\right).$$ 

Consider
\begin{eqnarray*}
Y_t^n 
&=& \chi_{E_n} \left(\mathbbm{E} [Y^{n-1} \, |\, \mathcal{F}^n_t] - Y^{n-1}_{T^n} + r_{T^n} X^{n-1}_{T^n}\right)  \\
&=& \chi_{E_n} \left( r_{T^n} X^{n-1}_{T^n} + \int^{t \vee T^n}_{T^n} \mathd Y^{n-1}_u \right) 
\end{eqnarray*}
and
\begin{eqnarray*}
 X_t^n &=& \chi_{E_n} \mathbbm{E} [X^{n-1} \; | \; \mathcal{F}^n_t]  \\
 &=& \chi_{E_n}  \left( X^{n-1}_{T^n} + \int^{t \vee T^n}_{T^n} \mathd X^{n-1}_u \right).
\end{eqnarray*}
By induction, these are martingales with respect to $\mathfrak{F}^n$ and $Y^n$ is
differentially subordinate to $X^n$.

\subsubsection{Sparseness} We prove that the resulting selection $\{T^ j\}_{j\geqslant 0}$ is sparse. Let $n\geqslant 0$. 
Let $A^{n}\subset E_{n}$ with $A^{n} \in \mathcal{F}_{T^{n}}$. Then $\chi_{A^{n}}Y^{n}_t $ is differentially 
subordinate to $\chi_{A^{n}}X^{n}_t $ in $\mathfrak{F}^{n}$.  Notice that $|X|_{T^{n}}=|X^{n}|_{0}$ since $\mathcal{F}^{n}_0=\mathcal{F}_{T^{n}}$.

Thanks to Lemma \ref{Theorem_WangWeakType} applied to the set $A=A_n$  there holds
$$\mathbbm{P} (A^{n} \cap E_{n+1})  = \{ \omega \in A^n : Y^{n \ast}(\omega) \vee X^{n \ast}(\omega) > 4 |X|_{T^n}(\omega)\} \leqslant  \frac{1}{2} \mathbb{P}(A^{n}).$$

\subsubsection{Domination} We prove the domination estimate. Indeed, we show that for all $n\geqslant 0$ there holds almost surely
\begin{equation}\label{estimate_inductive} Y^{\ast}\leqslant \sum_{j=0}^{n-1}8|X|_{T^j}\chi_{E_j}+Y^{n \ast}\quad (\mathcal{E}_n).\end{equation}
This implies the required domination because for the support of $Y^{n \ast}$ there holds
  $$\mathbb{P}\left(\text{supp} (Y^{n \ast})) \leqslant \mathbb{P}(E_{n}\right) \to 0 \text{ as } n\to \infty.$$
For $n=0$, the estimate $\mathcal{E}_0$ in (\ref{estimate_inductive}) follows from $$Y^{\ast}=Y^{\ast}\chi_{\Omega \setminus E_0} +  Y^{\ast}\chi_{E_0} = 0+ Y^{0 \ast}.$$
Assuming now $(\mathcal{E}_n)$ holds in (\ref{estimate_inductive}), we pass to $(\mathcal{E}_{n+1})$. Since $Y^{n \ast}$ is supported on $E_n$ we split 
$$Y^{n \ast} =  Y^{n \ast}\chi_{E_n \setminus E_{n+1}}+   Y^{n \ast}\chi_{E_{n+1}}.$$ 
In the complement of $E_{n+1}$ we have $Y^{n \ast} \leqslant 4 |X|_{T^n}\chi_{E_n}\leqslant 8 |X|_{T^n}\chi_{E_n}+Y^{n+1 \ast}$. In $E_{n+1}$ we have 
$$Y^{n \ast} (\omega)=\max \left\{ \sup_{t<T^{n+1}(\omega)} |Y^n_t(\omega)|, \sup_{t\geqslant T^{n+1}(\omega)} |Y^n_{t}(\omega)|   \right\}.$$
The first supremum is bounded by $4|X|_{T^n}\chi_{E_n}$ and for the second supremum we write trajectory-wise for $t\geqslant T^{n+1}(\omega)$
\begin{eqnarray*}
Y^n_t &=&  Y^n_0 + \int^{T^{n+1}}_{0} \mathd Y^n_u + \int^t_{T^{n+1}} \mathd Y^n_u  \\
 &=& Y^n_0 + \int^{T^{n+1}_-}_{0} \mathd Y^n_u 
+  (Y_{T^{n+1}}-Y_{T^{n+1}_-})  + \int^t_{T^{n+1}} \mathd Y^n_u  \\
&=&    \left( Y^n_0 + \int^{T^{n+1}_-}_{0} \mathd Y^n_u \right) - (r_{T^{n+1}} X^n_{T^{n+1}_-}) +
   \left( r_{T^{n+1}}  X^n_{T^{n+1}} + \int^t_{T^{n+1}} \mathd Y^n_u \right).
\end{eqnarray*}
We estimate in $E_{n+1}$ for $t\geqslant T^{n+1}$  
$$|Y^ n_t|  \leqslant  \left|Y^n_0 + \int^{T^{n+1}_-}_{0} \mathd Y^n_u \right| + | X^n_{T^{n+1}_-}| +  \left| r_{T^{n+1}} X_{T^{n+1}} + \int^t_{T^{n+1}} \mathd Y^n_u \right|.$$
      
The first two summands are each controlled by $4 | X |_{T^n}\chi_{E_n}$ by the definition of
the stopping time $T^{n+1}$. Last, observe that the third term on $E_{n+1}$ is
dominated by $ Y^{n+1 \ast}$.

Gathering the information, there holds almost surely
\[ Y^{n \ast} \leqslant 8 \, | X |_{T^n}\chi_{E_n} + Y^{n+1 \ast} , \] 
and thus
$$Y^{\ast} \leqslant \sum_{j=0}^{n-1} 8\,|X|_{T^j}\chi_{E_j} + Y^{n \ast} \leqslant \sum_{j=0}^{n} 8\,|X|_{T^j}\chi_{E_j} +Y^{n+1 \ast} .$$
The claim $(\mathcal{E}_{n+1})$ in (\ref{estimate_inductive}) is proved and the sparse domination in Theorem \ref{Theorem_sparsedomination} follows. 

\

\section{Maximal function of $Y$}

In this section we prove Theorems \ref{Theorem_SX} and \ref{Theorem_YstarX}. We first prove Theorem \ref{Theorem_SX} for $p=2$ and then obtain the result for other $p$ via extrapolation. We then deduce Theorem \ref{Theorem_YstarX} via the sparse domination Theorem \ref{Theorem_sparsedomination}. 

In order to prove Theorem \ref{Theorem_SX} 
for the case $p = 2$, it suffices to show that there exists $c_{\ref{Theorem_SX},2}\geqslant 0$ such that for all $w \in A_2$ and all functions $X \in L^2(w)$ there holds
\begin{equation} \label{sparse_weighted_estimate}
\| \mathcal{S} ( {X} ) \|_{L^2 (w)} \leqslant c_{\ref{Theorem_SX},2} Q_2 (w) \|
   {X} \|_{L^2 (w)} . 
\end{equation}   
This means 
\[ \left(\mathbbm{E} [(\mathcal{S} ( {X} ))^2 w]\right)^{\frac12} \leqslant c_{\ref{Theorem_SX},2} Q_2
   (w)\left( \mathbbm{E} [| {X} |^2 w] \right)^{\frac12}. \]
Dualizing and writing $u=w^{-1}$, we reduce to the estimate
\[ \mathbbm{E} [\mathcal{S} ({X} ) | {Z} |] \leqslant
   c_{\ref{Theorem_SX},2} Q_2 (w) \mathbbm{E} [| {X} |^2 w]^{\frac12} \mathbbm{E} [|{Z} |^2 u]^{\frac12}. \]
  We introduce the notations
\begin{equation}\label{weighted_expectation} 
\mathbbm{E} [\,\cdot \; w ]
  =\mathbbm{E}_w [\,\cdot \,] \mathbbm{E} [w ] .
   \end{equation}
  \begin{equation}\label{weighted_sampling} 
   \cdot_{\,\tau, w} =\mathbbm{E}_w [\, \cdot \;|\, 
\mathcal{F}_{\tau}].
  \end{equation}

Then, we write
$| \tilde{X} | u = | {X} |$ and $| \tilde{Z} | w = | {Z} |$ then suppressing the $\tilde \cdot$ again, it suffices
to prove
\[ \mathbbm{E} [\mathcal{S} (| X | u) | Z | w] 
\leqslant c_{\ref{Theorem_SX},2} Q_2 (w)
   \mathbbm{E} [w]^{\frac12} \mathbbm{E} [u]^{\frac12} 
   \mathbbm{E}_u [| X |^2]^{\frac12} \mathbbm{E}_w [| Z |^2]^{\frac12} . \]
Now, we calculate the left hand side
\begin{eqnarray*}
 \mathbbm{E} \left[ \sum_j (| X | u)_{T^j} \chi_{E_j} | Z | w \right]
   &=&\mathbbm{E} \left[ \sum_j \mathbbm{E} [(| X | u)_{T^j} \chi_{E_j} | Z | w \;
   | \; \mathcal{F}_{T^j}] \right] \\
   &=&\mathbbm{E} \left[ \sum_j (| X | u)_{T^j} (|
   Z | w)_{T^j} \chi_{E_j} \right] \\
  &\leqslant& Q_2 (w) \mathbbm{E} \left[ \sum_j | X |_{T^j, u} | Z |_{T^j, w} \chi_{E_j}
   \right]. 
 \end{eqnarray*}  
In the above calculation, we used the notations (\ref{weighted_expectation}) and (\ref{weighted_sampling}) and noticed that
\begin{equation*}
(| Z | w)_{T^j} \chi_{E_j} =\mathbbm{E} [| Z | w\, |\, \mathcal{F}_{T^j}]
   \chi_{E_j} =\mathbbm{E}_w [| Z | | \mathcal{F}_{T^j}] \mathbbm{E} [w \,|\,
   \mathcal{F}_{T^j}] \chi_{E_j} ,
   \end{equation*}
and similarly for the other term. We recalled that by the $A_2$ condition
\[ \| \mathbbm{E} [w \,|\, \mathcal{F}_{T^j}] \mathbbm{E} [w^{- 1} \,|\,
   \mathcal{F}_{T^j}] \|_{\infty} \leqslant Q_2 (w). \]
For each fixed $j$ we have that the non-negative random variable $$| X
|_{T^j, u} | Z |_{T^j, w} \chi_{E_j} \in \mathcal{F}_{T^j}$$ and as such it can
be approximated from below by step functions.
\[ \sum_k \alpha_k^j \chi_{A_k^j} \nearrow | X |_{T^j, u} | Z |_{T^j, w}
   \chi_{E_j} ,\]
with $A^j_k \in \mathcal{F}_{T^j}$ disjoint and $\stackrel{\cdot}{\cup}_k A^j_k = E_j$. Notice
that on $A^j_k$ there holds
\begin{equation}\label{estimate_maximal}
 \alpha^j_k \chi_{A^j_k}(\omega)\leqslant X_{u}^{\ast} Z_{ w}^{\ast} (\omega),
 \end{equation}
where the maximal functions are taken with respect to 
 weighted measure:
\begin{equation}\label{notation_weightedmax} 
X_{ u}^{\ast}(\omega)=\sup_{t}| X_{t,u}(\omega)|.
\end{equation}

\

Now recall that $\mathbbm{P} (A_k^j \cap E_{j + 1}) \leqslant \frac{1}{2}
\mathbbm{P} (A_k^j)$ and so if we write  $$S^{j+1}_k = A_k^j \backslash (A_k^j \cap E_{j + 1}),$$ then
$\mathbbm{P} (S_k^{j+1} ) \geqslant \frac{1}{2} \mathbbm{P} (A_k^j)$. 
We changed the index on $S$ to recall the important fact that it is measurable in $\mathcal{F}_{T^{j+1}}$. 
Notice the crucial property of the collection $\{S_k^{j+1}\}_{j,k\geqslant 0}$ : it is a disjoint collection in both parameters.

\

We estimate 
\begin{eqnarray*}
  &&\mathbbm{E} \left[ \sum^J_{j = 0} \sum_k \alpha_k^j \chi_{A_k^j} \right] 
  =  
  \sum^J_{j = 0} \sum_k \alpha_k^j \mathbbm{E} \left[\chi_{A^j_k}\right]
   =  \sum^J_{j = 0} \sum_k \alpha_k^j \mathbbm{P} \left(A^j_k\right)\\
  & \leqslant & 
  2 \sum^J_{j = 0} \sum_k \alpha^j_k \mathbbm{P} \left(S^{j+1}_k\right)
   =  2\,\mathbbm{E} \left[ \sum^J_{j = 0} \sum_k \alpha_k^j \chi_{S^{j+1}_k} \right]\\
  & = & 
  2\,\mathbbm{E} \left[ \sum^J_{j = 0} \sum_k \alpha_k^j\chi_{S^{j+1}_k} w^{\frac12} u^{\frac12} \right]
  \leqslant  
  \mathbbm{E} \left[ \sum^J_{j = 0} \sum_k X^{\ast}_{ u} u^{\frac12} Z^{\ast}_{ w} w^{\frac12} \chi_{S^{j+1}_k}  \right]\\
   &\leqslant&  
   2 \left( \mathbbm{E} \left[ \sum^J_{j = 0} \sum_k (X^{\ast}_{ u})^2 u \chi_{S^{j+1}_k} \right] \right)^{\frac12} 
    \left(\mathbbm{E} \left[ \sum^J_{j = 0} \sum_k (Z^{\ast}_{ w})^2 w \chi_{S^{j+1}_k} \right] \right)^{\frac12}\\[.5em] 
  & \leqslant & 
  2 \left(\mathbbm{E} \left[\left(X^{\ast}_u\right)^2 u\right]\right)^{\frac12} \left(\mathbbm{E}
  \left[\left(Z_w^{\ast}\right)^2 w\right]\right)^{\frac12}\\[.5em] 
  & = & 
  2 \left(\mathbbm{E} [u]\right)^{\frac12} \left(\mathbbm{E} [w]\right)^{\frac12} \left(\mathbbm{E}_u
  \left[\left(X^{\ast}_u\right)^2\right]\right)^{\frac12} \left(\mathbbm{E}_w \left[\left(Z^{\ast}_w\right)^2\right]\right)^{\frac12}\\[.5em] 
  & \leqslant & 
  8 \,(\mathbbm{E} [w])^{\frac12} (\mathbbm{E} [u])^{\frac12}
  \left(\mathbbm{E}_u \left[| X |^2\right]\right)^{\frac12} \left(\mathbbm{E}_w \left[| Z |^2\right]\right)^{\frac12}.
\end{eqnarray*}
By the monotone convergence theorem, this gives us the estimate
\[ \mathbbm{E} \left[ \sum_j | X |_{T^j, u} | Z |_{T^j, w} \chi_{E_j} \right]
   \leqslant 8\,\mathbbm{E} [w]^{\frac12} \mathbbm{E} [u]^{\frac12} \mathbbm{E}_u \left[|
   X |^2\right]^{\frac12} \mathbbm{E}_w \left[| Z |^2\right]^{\frac12}, \]
   and  we have thus seen that inequality (\ref{sparse_weighted_estimate}) holds with $c_{\ref{Theorem_SX},2}=8$.

   \
   
We now point out the changes for the case $p \neq 2$. We use the extrapolation theorem below from \cite{DomPetWit2017a}. 
Notice that only the weights are required to have the martingale property while $X,Y$ are functions.

\begin{theorem}[Domelevo-Petermichl]\label{thmextrap}
  Given a filtered probability space as described above. Let $1 < p < \infty$
  and $w \in A_p$. Let  $X, Y \in L^p(w)$. Suppose
  $1 < r < \infty$ and suppose $\forall A \geqslant 1 \; \exists N_r (A) > 0$
  increasing such that for triples $X, Y, \rho$ with $X,Y \in L^r(\rho)$ and
  $Q_r (\rho) \leqslant A$ there holds
  \begin{eqnarray*}
    \| Y \|_{L^r(\rho)} \leqslant N_r (A) \| X \|_{L^r(\rho)} .
  \end{eqnarray*}
  Then for any $1 < p < \infty$ there exists $N_p (B) > 0$ such that if
  $Q_p (w)  \leqslant B$ there holds
  \begin{eqnarray*}
    \| Y \|_{L^p(w)} \leqslant N_p (B) \| X \|_{L^p(w)} .
  \end{eqnarray*}
  With $c_{\ref{Theorem_XstarX},p}$ denoting the numeric part of the estimate in the
  weighted $L^p$ maximal estimate  from Theorem \ref{Theorem_XstarX}, in particular 
   $$N_p (B) \leqslant 2^{\frac1{r}} N_r \left(2 c_{\ref{Theorem_XstarX},p'}^{\frac{p - r}{p - 1}} B\right) \text{ if } p > r .$$ 
  $$N_p (B) \leqslant 2^{\frac{r - 1}{r}} N_r \left( 2^{r - 1} \left(c_{\ref{Theorem_XstarX},p}^{p - r} B\right)^{\frac{r - 1}{p - 1}} \right) \text{ if } p < r. $$
\end{theorem}

Using this theorem for $r=2$ we extrapolate from inequality (\ref{sparse_weighted_estimate}) 
\[ \| \mathcal{S} ( {X} ) \|_{L^2 (w)} \leqslant 8\, Q_2 (w) \| {X} \|_{L^2 (w)} , \]
so $N_2(A)=8A$. We obtain  the estimate claimed in Theorem \ref{Theorem_SX}
$$\| \mathcal{S} ( X ) \|_{L^p (w)} \leqslant c_{\ref{Theorem_SX},p} Q_p (w)^{\max \{1, \frac1{p-1}\}} \| X \|_{L^p (w)}.$$ 
Finally 
$$\| Y^{\ast} \|_{L^p (w)} \leqslant 8\| \mathcal{S} ( X ) \|_{L^p (w)}$$ 
gives the claimed estimate in Theorem \ref{Theorem_YstarX}.

\section{Maximal Function of $X$}

In this section we prove Theorem \ref{Theorem_XstarX} via the modification of a simple and direct
 domination argument for the maximal function. See for example the argument by Lerner {\cite{Ler2008a}} in a different context. Notice that the obtained norm estimate is the same as that in 
 Buckley's text \cite{Buc1993a} on homogeneous spaces. Buckley's proof enjoyed an extension to some martingales with certain restrictive homogeneity 
 conditions in the presence of jumps - this was needed because of the failure of the openness of the $A_p$ class in the general context. The argument here 
 does not rely on the openness condition of the $A_p$ classes and is therefore providing the estimate in full generality and in addition recovers the correct growth with the $A_p$ characteristic for the norm estimate. The argument consists of a trajectory-wise domination of the maximal operator of $X$ and the use of Doob's inequality. We write $u=w^{\frac{1}{p-1}}$ for the dual weight and recall that the $A_p$ characteristic is $\sup_{\tau}\| u_{\tau}^{p-1}w_{\tau}\|_{\infty}$. We remind the reader of notations (\ref{weighted_expectation}), (\ref{weighted_sampling}) and (\ref{notation_weightedmax}). There holds for all $t\ge0$ and all $p>1$ 
\begin{eqnarray*}
  | X_t |^{p - 1} 
  & \leqslant &
   \left(\mathbbm{E} \left[| X |\, | \, \mathcal{F}_t\right]\right)^{p -1}\ =  
  \left(\mathbbm{E} \left[| X | u^{-1}u \, | \, \mathcal{F}_t\right]\right)^{p - 1}\\
  & = & 
  \left(\mathbbm{E}_u \left[| X | u^{-1} \, | \, \mathcal{F}_t\right]\right)^{p - 1} \left(\mathbbm{E}
  \left[u \, | \, \mathcal{F}_t\right]\right)^{p - 1}\\
  & \leqslant & 
  Q_p (w) \left(\mathbbm{E} \left[w  \, | \, \mathcal{F}_t\right]\right)^{- 1}
  \left(\mathbbm{E}_u \left[| X | u^{-1} \, | \,  \mathcal{F}_t\right]\right)^{p - 1}.
\end{eqnarray*}
Now observe that 
\begin{eqnarray*}
\left(\mathbbm{E}_u \left[| X | u^{-1}\, |\, \mathcal{F}_t\right]\right)^{p - 1}
&=&
\mathbbm{E} \left[\left(\mathbbm{E}_u [| X | u^{-1} \,  | \, \mathcal{F}_t]\right)^{p - 1} | \,\mathcal{F}_t\right] \\
&\leqslant& 
\mathbbm{E} \left[\left(\left(| X | u^{-1}\right)_{^{} u}^{\ast}\right)^{p - 1}w^{-1} w \, | \, \mathcal{F}_t\right]\\
&=&
\mathbbm{E}_w\left[ \left(\left(| X | w\right)_{^{} u}^{\ast}\right)^{p - 1}w^{-1}  \, | \, \mathcal{F}_t\right]  
\mathbbm{E}\left[ w \, | \, \mathcal{F}_t \right] .
\end{eqnarray*} 
Then get for all $t$
\[ | X_t |^{p - 1} \leqslant Q_p (w) \left(\left(\left(| X | u^{- 1}\right)^{\ast}_u\right)^{p - 1} w^{-1} \right)^{\ast}_w,\]
and therefore 
$$(X^{\ast })^p\leqslant \left( Q_p (w) \left(\left(\left(| X | u^{- 1}\right)^{\ast}_u\right)^{p - 1} w^{-1} \right)^{\ast}_w \right)^{\frac{p}{p-1}}.$$
Thus
\begin{eqnarray*}
  \mathbbm{E} [(X^{\ast})^p w] 
  & \leqslant & Q_p (w)^{\frac{p} {p - 1}} 
  \mathbbm{E}\left[\left(\left(\left(\left(| X | u^{- 1}\right)^{\ast}_u\right)^{p - 1} w^{- 1}\right)^{\ast}_w\right)^{p'} w\right]\\
  & = & Q_p (w)^{\frac{p} {p - 1}} 
  \mathbbm{E}_w \left[\left(\left(\left(\left(| X | u^{- 1}\right)^{\ast}_u\right)^{p -1} w^{- 1}\right)^{\ast}_w\right)^{p'}\right] \mathbbm{E} [w]\\
  & \leqslant & Q_p (w)^{\frac{p} {p - 1}}  \left( \frac{p'}{p' - 1} \right)^{p'}
  \mathbbm{E}_w \left[(((| X | u^{- 1})^{\ast}_u)^{p - 1} w^{- 1})^{p'}\right]
  \mathbbm{E} [w]\\
  & = & Q_p (w)^{\frac{p} {p - 1}}  \left( \frac{p'}{p' - 1} \right)^{p'} \mathbbm{E}
  \left[(((| X | u^{- 1})^{\ast}_u)^{p - 1} w^{- 1})^{p'} w\right]\\
  & = & Q_p (w)^{\frac{p} {p - 1}} \left( \frac{p'}{p' - 1} \right)^{p'} \mathbbm{E}
 \left[((| X | u^{- 1})^{\ast}_u)^p u\right]\\
  & = & Q_p (w)^{\frac{p} {p - 1}}  \left( \frac{p'}{p' - 1} \right)^{p'}
  \mathbbm{E}_u \left[((| X | u^{- 1})^{\ast}_u)^p\right] \mathbbm{E} [u]\\
  & \leqslant & Q_p (w)^{\frac{p} {p - 1}} \left( \frac{p'}{p' - 1} \right)^{p'}
  \left( \frac{p}{p - 1} \right)^p \mathbbm{E}_u \left[(| X | u^{- 1})^p\right]
  \mathbbm{E} [u]\\
  & = & Q_p (w)^{\frac{p} {p - 1}}  \left( \frac{p'}{p' - 1} \right)^{p'} \left(
  \frac{p}{p - 1} \right)^p \mathbbm{E} \left[(| X | u^{- 1})^p u\right]\\
  & = & Q_p (w)^{\frac{p} {p - 1}}  \left( \frac{p'}{p' - 1} \right)^{p'} \left(
  \frac{p}{p - 1} \right)^p \mathbbm{E} \left[| X |^p w\right].
\end{eqnarray*}
Raising to the power $1/p$ gives the desired estimate in Theorem \ref{Theorem_XstarX} with $c_{\ref{Theorem_XstarX},p}=\frac{p^{p'}}{p-1}$.

\


\begin{thebibliography}{10}
  \bibitem{BekBon1978a}David Bekoll{\'e} and Aline Bonami.
  {\newblock}In{\'e}galit{\'e}s {\`a} poids pour le noyau de Bergman.
  {\newblock}\tmtextit{C. R. Acad. Sci. Paris S{\'e}r. A-B},
  286(18): A775--A778, 1978.
  
  \bibitem{BonLep1979a}Aline Bonami and Dominique L\'{e}pingle. {\newblock}Fonction
  maximale et variation quadratique des martingales en pr{\'e}sence d'un
  poids. {\newblock}\tmtextit{S{\'e}minaire de probabilit{\'e}s XIII, Univ. Strasbourg
  1977/78, Lect. Notes Math.} 721: 294--306, 1979.
  
  \bibitem{Buc1993a}Stephen~M. Buckley. {\newblock}Estimates for operator
  norms on weighted spaces and reverse Jensen inequalities.
  {\newblock}\tmtextit{Trans. Amer. Math. Soc.}, 340(1): 253--272, 1993.
  
   \bibitem{Bur1991a}Donald~L. Burkholder.
\newblock Explorations in martingale theory and its applications.
{\newblock}\tmtextit{\'{E}cole d'\'{E}t{\'e} de Probabilit{\'e}s de
  Saint-Flour}, XIX, Lecture Notes in Math.,1464:  
  1--66, 1991.

  
  \bibitem{CruMarPer2012a}David Cruz-Uribe, Jos{\'e} Martell, and Carlos
  P{\'e}rez. {\newblock}Sharp weighted estimates for classical operators.
  {\newblock}\tmtextit{Adv. Math.}, 229(1): 408--441, 2012.
  
   \bibitem{DahDomPetSkr2019}Kamilia Dahmani, Komla Domelevo, 
   Stefanie Petermichl, and Kristina \v{S}kreb. {\newblock} 
   Dimensionless weighted estimates for the Bakry Riesz vector. {\newblock} in preparation 2019.
  
  \bibitem{DelMey1982a}Claude Dellacherie and Paul-Andr{\'e} Meyer.
  {\newblock}Probabilities and potential. B, volume~72 of
  \tmtextit{North-Holland Mathematics Studies}. {\newblock}North-Holland
  Publishing Co., Amsterdam, 1982. {\newblock}Theory of martingales,
  Translated from the French by J. P. Wilson.
  
  \bibitem{DIPTV2019}Komla Domelevo, Paata Ivanisvili, Stefanie Petermichl, Sergei Treil, and Alexander Volberg.
  {\newblock}On the failure of lower square function estimates in the
  non-homogeneous weighted setting.
  {\newblock}\tmtextit{Math. Ann.}{\newblock}https://doi.org/10.1007/s00208-018-1787-4.
   2019.
  
  \bibitem{DolMey1979a}Catherine Dol{\'e}ans-Dade and Paul-Andr\'{e} Meyer.
  {\newblock}In{\'e}galit{\'e}s de normes avec poids. {\newblock}In
  \tmtextit{S{\'e}minaire de Probabilit{\'e}s, XIII (Univ. Strasbourg,
  Strasbourg, 1977/78)}, volume 721 of \tmtextit{Lecture Notes in Math.}: 
  313--331. Springer, Berlin, 1979.
  
   \bibitem{DomPetWit2017a}Komla Domelevo, Stefanie Petermichl, and Janine Wittwer.
  {\newblock}A dimensionless weighted bound for the Riesz vector in $\mathbb{R}^n$.
  {\newblock}\tmtextit{Bull. Sci. Math.}, 141(5): 385--407, 2017.
  
  \bibitem{DomPet2016a}Komla Domelevo and Stefanie Petermichl.
  {\newblock}Differential subordination under change of law.
  {\newblock}\tmtextit{Ann. Prob.}, 47(2): 896--925, 2019.
  
  \bibitem{GunVar1979}Richard~F. Gundy and Nicolas~Th. Varopoulos.
  {\newblock}Les transformations de Riesz et les int{\'e}grales stochastiques.
  {\newblock}\tmtextit{C. R. Acad. Sci. Paris S{\'e}r. A-B}, 289(1): A13--A16,
  1979.
  
  \bibitem{Hyt2012a}Tuomas~P. Hyt{\"o}nen. {\newblock}The sharp weighted
  bound for general Calder{\'o}n-Zygmund operators. {\newblock}\tmtextit{Ann.
  Math. (2)}, 175(3): 1473--1506, 2012.
  
  \bibitem{IzuKaz1977a}Masataka Izumisawa and Norihiko Kazamaki. {\newblock}Weighted
  norm inequalities for martingales. {\newblock}\tmtextit{T{\^o}hoku Math. J.
  (2)}, 29(1): 115--124, 1977.
  
  \bibitem{Lac2015a}Michael~T. Lacey. {\newblock}An elementary proof of
  the $A_2$ bound. {\newblock}\tmtextit{Israel Journal of Mathematics}, 217: 181-195, 2017.
  
  \bibitem{LacPetReg2010a}Michael~T. Lacey, Stefanie Petermichl, and
  Maria~C. Reguera. {\newblock}Sharp $A_2$ inequality for Haar shift
  operators. {\newblock}\tmtextit{Math. Ann.}, 348(1): 127--141, 2010.
    
  \bibitem{Ler2008a}Andrei~K. Lerner. {\newblock}An elementary approach to
  several results on the Hardy-Littlewood maximal operator.
  {\newblock}\tmtextit{Proc. Amer. Math. Soc.}, 136(8): 2829--2833, 2008.
  
  \bibitem{LerNaz}Andrei~K. Lerner and Fedor Nazarov.{\newblock}
  Intuitive dyadic calculus.{\newblock}\tmtextit{ArXiv},1508.05639: 1--53, 2015.
  
  \bibitem{Muc1972a}Benjamin Muckenhoupt. {\newblock}Weighted norm
  inequalities for the Hardy maximal function. {\newblock}\tmtextit{Trans.
  Amer. Math. Soc.}, 165: 207--226, 1972.
  
  \bibitem{NazTreVol1999a}Fedor Nazarov, Sergei Treil, and Alexander Volberg.
  {\newblock}The Bellman functions and two-weight inequalities for Haar
  multipliers. {\newblock}\tmtextit{J. Amer. Math. Soc.}, 12(4): 909--928,
  1999.
  
  \bibitem{Ose2016a}Adam Os\c{e}kowski. {\newblock}Sharp $L^p$-bounds for the
  martingale maximal function. {\newblock}\tmtextit{to appear in Tohoku
  Mathematical Journal}, 2016.
  
  \bibitem{Pet2000b}Stefanie Petermichl. {\newblock}Dyadic shifts and a
  logarithmic estimate for Hankel operators with matrix symbol.
  {\newblock}\tmtextit{C. R. Acad. Sci. Paris S{\'e}r. I Math.},
  330(6): 455--460, 2000.
  
  \bibitem{PetVol2002a}Stefanie Petermichl and Alexander Volberg.
  {\newblock}Heating of the Ahlfors-Beurling operator: weakly quasiregular
  maps on the plane are quasiregular. {\newblock}\tmtextit{Duke Math. J.},
  112(2): 281--305, 2002.
  
  \bibitem{Pro2005a}Philip~E. Protter. {\newblock}\tmtextit{Stochastic
  integration and differential equations}, volume~21 of \tmtextit{Stochastic
  Modelling and Applied Probability}. {\newblock}Springer-Verlag, Berlin,
  2005. {\newblock}Second edition. Version 2.1, Corrected third printing.
  
  \bibitem{ThiTreVol2015a}Christoph Thiele, Sergei Treil, and Alexander
  Volberg. {\newblock}Weighted martingale multipliers in the non-homogeneous
  setting and outer measure spaces. {\newblock}\tmtextit{Adv. Math.},
  285: 1155--1188, 2015.
  
  \bibitem{Wan1995a}Gang Wang. {\newblock}Differential subordination and
  strong differential subordination for continuous-time martingales and
  related sharp inequalities. {\newblock}\tmtextit{Ann. Probab.},
  23(2): 522--551, 1995.
  
  \bibitem{Wit2000a}Janine Wittwer. {\newblock}A sharp estimate on the
  norm of the martingale transform. {\newblock}\tmtextit{Math. Res. Lett.},
  7(1): 1--12, 2000.
\end{thebibliography}
\end{document}